\documentstyle[12pt]{article}
\input amssym.def

\newfont{\fvar}{msbm10 at 14pt}

\def\0{\hbox{\bf 0}}
\def\1{\hbox{\bf 1}}
\def\bR{\hbox{\bf R}}

\begin{document}

Doklady Mathematics, vol. 58, No. 3, 1998, pp. 389--391

Translated from Doklady Akademii Nauk,

vol. 363, No. 3, 1998, pp. 298--300

\bigskip

\begin{center}
{\Large
Linear Functionals on Idempotent Spaces:

An Algebraic Approach}
\footnote{International Sophus Lie Centre,
Moscow, Russia, e-mail:islc@islc.msk.su\qquad}

\medskip

G.L. Litvinov, V.~P.~Maslov and  G.~B.~Shpiz
\end{center}

\bigskip

In this paper, we present an algebraic approach to idempotent
functional analysis, which is an abstract version of idempotent
analysis in the sense of [1--3]. Elements of such an
approach were used, for example, in [1, 4]. The basic concepts and
results are expressed in purely algebraic terms. We consider
idempotent versions of certain basic results of linear
functional analysis, including the theorem on the general
form of a linear functional and the Hahn--Banach and
Riesz--Fischer theorems.

1. Recall that an additive semigroup $S$ with commutative addition
$\oplus$ is called an {\it idempotent semigroup} (IS) if the relation
$x\oplus x=x$ is fulfilled for all elements $x\in S$. If $S$
contains a neutral element, this element is denoted by the
symbol $\0$. Any IS is a partially ordered set with respect
to the following standard order: $x\preceq y$ if and only
if $x\oplus y=y$. It is obvious that this order is well defined
and $x\oplus y=\sup \{x, y\}$. Thus, any IS is an upper semilattice;
moreover, the concepts of IS and upper semilattice coincide
[5]. An idempotent semigroup $S$ is called $a$-{\it complete}
(or {\it algebraically complete}) if it is complete as an ordered set,
i.e., if any subset $X$ in $S$  has the least upper bound
$\sup(X)$ denoted by $\oplus X$ and the greatest lower bound $\inf (X)$
denoted by $\wedge X$. This semigroup is called $b$-{\it complete}
(or {\it boundedly complete}), if any bounded above subset $X$ of this
semigroup (including the empty subset) has the least upper bound
$\oplus X$ (in this case, any nonempty subset $Y$ in $S$
has the greatest lower bound $\wedge Y$ and $S$ in a lattice). Note that
any $a$-complete or $b$-complete IS has the zero element $\0$
that coincides with $\oplus{\fvar\symbol{"1F}}$,
where ${\fvar \symbol{"1F}}$ is the empty set. Certainly,
$a$-completeness implies the $b$-completeness.
Completion by means of cuts [5] yields an embedding
$S\to \widehat S$ of an arbitrary IS $S$ into an $a$-complete IS
$\widehat S$ (which is called a {\it normal completion} of $S$); in addition,
$\widehat{\widehat S} = S$. The $b$-completion procedure
$S \to \widehat S_b$ is defined similarly: if $S \ni \infty =\sup S$,
then $\widehat S_b =S$; otherwise,
$\widehat S =\widehat S_b \cup \{ \infty \}$.
An arbitrary $b$-complete IS  $S$ also may differ from $\widehat S$ only
by the element $\infty =\sup S$.

Let $S$ and $T$ be $b$-complete IS.
Then, a homomorphism $f : S\to T$ is said to be a $b$-{\it homomorphism} if
 $f (\oplus X) = \oplus f(X)$ for any bounded subset $X$ in $S$.
If the  $b$-homomorphism $f$ is extended to a homomorphism
$\widehat S\to \widehat T$ of the correesponding normal completions and
$f(\oplus X) = \oplus f(X)$ for all $X\subset S$, then $f$
is said to be an $a$-{\it homomrphism}. An IS $S$ equipped with a topology such
that the set $\{ s\in S\vert s\preceq b\}$ is closed in this topology
for any $b\in S$ is called a {\it topological idempotent semigroup} $S$.

{\bf Proposition 1.} {\sl Let $S$ be an $a$-complete topological
IS and $T$ be a  $b$-complete topological IS such that, for
any nonempty subsemigroup $X$ in $T$, the element $\oplus X$ is
contained in the topological closure of $X$ in $T$. Then, a
homomorphism $f : T\to S$ that maps zero into zero
is an $a$-homomorphism if and only
if the mapping $f$ is lower semicontinuous in the sense that
the set $\{ t\in T\vert f(t)\preceq s\}$ is closed in $T$ for any} $s\in S$.

{\bf 2.} An idempotent semigroup $K$ on which an associative
multiplication $\odot$ with identity $\1$ is defined together
with the idempotent addition $\oplus$ and both distributive relations
are fulfilled is called an {\it idempotent semiring} (ISR). The element
$\0\in K$, $\0\neq\1$ is called a {\it zero} of the semiring $K$  if
$x\oplus\0=x$ and $x\odot\0=\0\odot x=\0$ for any $x\in K$. A
commutative ISR in which every nonzero element is invertible
with respect to the multiplication is called an {\it idempotent
semifield} (or briefly, {\it semifield}). An idempotent semiring
$K$ is called $a$-{\it complete} (respectively $b$-{\it complete})
if $K$ is an $a$-complete (respectively $b$-complete) IS and,
for any subset (respectively, for any bounded subset) $X$ in $K$
and any $k\in K$, the generalized distributive laws
$k\odot(\oplus X)=\oplus (k\odot X)$ and $(\oplus X)\odot k =
\oplus(X\odot k)$ are fulfilled. Generalized distributivity
implies that any $a$-complete or $b$-complete ISR has a zero
element that coincides with $\oplus{\fvar \symbol{"1F}}$,
where ${\fvar \symbol{"1F}}$ is the empty set.
The concept of $a$-complete ISR coincides with the concept of
complete dioid in the sense of [4].

The set $\bR(\max, +)$ of real numbers equipped with the idempotent
addition $\oplus=\max$ and multiplication $\odot=+$ is an example of an
ISR; in this case, $\1 = 0$. Adding the element $\0=-\infty$
to this ISR, we obtain a $b$-complete semiring
$\bR_{\max} = \bR\cup \{-\infty\}$
with the same operations and the same zero.
Adding the element $=\infty$ to $\bR_{\max}$ and asumming that
$\0\odot (+\infty)=\0$ and $x\odot (+\infty) = +\infty$
for $x\neq\0$ and $x\oplus (+\infty) = +\infty$ for any $x$, we
obtain the $a$-complete ISR $\widehat{\bR}_{\max} = \bR_{\max}\cup \{+\infty\}$.
The standard order on $\bR(\max, +)$, $\bR_{\max}$ and
$\widehat{\bR}_{\max}$ coincides with the ordinary order.
The ISRs $\bR(\max,+)$ and  $\bR_{\max}$ are semifields.
On the contrary, an $a$-complete ISR that does not coincide with
$\{\0, \1\}$ cannot be a semifield. An important class
of examples is related to
(topological) vector lattices (see, for example, [5];
[6, Chapter 5]). Defining the sum $x\oplus y$ as $\sup\{x, y\}$ and
the multiplication $\odot$ as the addition of vectors, we can
interpret the vector lattices as idempotent semifields.
Adding the element $\0$ to a complete vector lattice
(in the sense of [5,6]), we obtain a $b$-complete semifield.
If, in addition, we add the infinite element, we obtain an
$a$-complete ISR (which, as an ordered set, coincides with the normal
completion of the original lattice).

{\bf 3.} Let $V$ be an idempotent semigroup and  $K$ be an
idempotent semiring. Suppose that a multiplication $k, x\mapsto k\odot x$
of all elements from $K$ by the elements from
$V$ is defined; moreover, this multiplication is associative
and distributive with respect
to the addition in $V$ and $\1\odot x = x$, $\0\odot x=\0$
for all $x\in V$. In this case, the semigroup $V$ is called an
{\it idempotent semimodule} (or simply, a {\it semimodule})
 over $K$. The element
 $\0_V\in V$ is called the {\it zero} of the semimodule $V$ if $k\odot\0_V=\0_V$
and $\0_V\oplus x = x$ for any $k\in K$ and $x\in V$.
Let $V$ be a semimodule over a $b$-complete idempotent semiring $K$.
This semimodule is called $b$-{\it complete} if it is $b$-complete as
an IS and, for any bounded subsets $Q$ in $K$ and $X$ in  $V$,
the generalized distributive laws $(\oplus Q)\odot x = \oplus (Q\odot x)$ and
$k\odot (\oplus X) = \oplus (k\odot X)$ are fulfilled for all
$k\in K$ and $x\in X$. This semimodule is called $a$-{\it complete}
if it is $b$-complete and contains the element $\infty = \sup V$.

A semimodule $V$ over a $b$-complete semifield $K$ is said to be
an $a$-{\it idempotent} $a$-{\it space} ($b$-{\it space})
 if this semimodule is
$a$-complete (respectively, $b$-complete) and the equality
$(\wedge Q)\odot x = \wedge (Q\odot x)$ holds for any nonempty subset
$Q$ in $K$ and any $x\in V$, $x\neq \infty = \sup V$.
The normal completion $\widehat V$ of a $b$-space (as an IS) has the
structure of an idempotent $a$-space (and may differ from $V$ only
by the element $\infty= \sup V$).

Let $V$ and $W$ be idempotent semimodules over an idempotent semiring
$K$.  A mapping $p: V\to W$ is said to be {\it linear} (over $K$) if
$$
p(x\oplus y)=p(x)\oplus p(y) \mbox{ and } p(k\odot x)=k\odot p(x)
$$
for any $x, y\in V$ and $k\in K$. Let the semimodules $V$ and $W$ be
 $b$-complete. A linear mapping $p: V\to W$ is said to be $b$-{\it linear}
if it is a $b$-homomorphism of the IS; this mapping is said to be
$a$-{\it linear} if it can be extended to an $a$-homomorphism of the normal
completions $\widehat V$ and $\widehat W$. Proposition 1 (see above)
shows that $a$-linearity simulates (semi)continuity for
linear mappings. The normal completion $\widehat K$ of the
semifield $K$ is a semimodule over $K$. If $W= \widehat K$, then the
linear mapping $p$ is called a linear functional.

Examples of idempotent semimodules and spaces that are the most important
for analysis are either subsemimodules of topological vector lattices
[6] (or coincide with them) or are dual to them, i.e., consist
of linear functionals subject to some regularity condition,
for example, consist of $a$-linear functionals.

{\bf 4.} Let $V$ be an idempotent $b$-space over  a
$b$-complete semifield $K$,  $x\in \widehat V$. Denote by
 $x^*$ the functional $V\to \widehat K$ defined by the formula
$x^* (y)=\wedge \{ k\in K \vert y\preceq k\odot x\}$,
where $y$ is an arbitrary fixed element from $V$.

{\bf Theorem 1.} {\sl For any $x\in \widehat V$ the functional $x^*$
is $a$-linear. Any nonzero $a$-linear functional $f$ on $V$
is given by $f = x^*$ for a unique suitable element $x\in V$.
If $K\neq \{\0, \1\}$, then} $x=\oplus\{ y\in V\vert f(y)\preceq\1\}$.

Note that results of this type obtained earlier concerning the structure of
linear functionals cannot be carried over to subspaces and
subsemimodules.

A subsemigroup $W$ in $V$ closed with respect to the
multiplication by an arbitrary element from $K$ is called a
$b$-{\it subspace} in $V$ if the imbedding $W\to V$ can be extended to a
$b$-linear mapping. The following result is obtained from Theorem 1 and
is the idempotent analog of the Hahn--Banach theorem.

{\bf Theorem 2.} {\sl Any $a$-linear functional defined on a
$b$-subspace $W$ in $V$ can be extended to an $a$-linear functional
on $V$. If $x, y\in V$ and $x\neq y$, then there exists an 
$a$-linear functional $f$ on $V$ that separates the elements
$x$ and $y$, i.e.,} $f(x)\neq f(y)$.

The following statements are easily derived from the definitions and 
can be regarded as the analogs of the well-known results of ordinary
functional analysis (the Banach--Steinhaus and the closed-graph
theorems).

{\bf Proposition 2.} {\sl Suppose that $P$ is a family of $a$-linear
mappings of an $a$-space $V$ into an $a$-space $W$ and the mapping
$p : V\to W$ is the pointwise sum of the mappings of this family,
i.e., $p(x) =\sup \{ p_{\alpha} (x)\vert p_{\alpha} \in P\}$. 
Then the mapping $p$ is $a$-linear.}

{\bf Proposition 3.} {\sl Let $V$ and $W$ be  $a$-spaces. A linear
mapping $p : V \to W$ is $a$-linear if and only if its graph
$\Gamma$ in $V\times W$ is closed with respect to passing
to sums {\rm (}i.e., to  least upper bounds{\rm )}
of its arbitrary subsets.}

{\bf 5.} Let $K$ be a $b$-complete semifield and $A$ be an
idempotent $b$-space over $K$ equipped with the structure of a semiring
compatible with the multiplication $K\times A\to A$ so that
the associativity of the multiplication is preserved. In this case,
$A$ is called an {\it idempotent $b$-semialgebra} over~$K$.

{\bf Proposition 4.} {\sl For any invertible element $x\in A$ from the
$b$-semialgebra $A$ and any element $y\in A$, the equality
$x^* (y) = \1^* (y\odot x^{-1})$ holds, where} $\1\in A$.

The mapping $A\times A\to \widehat K$ defined by the formula
$(x, y) \mapsto \langle x, y\rangle = \1^* (x\odot y)$
is called the  {\it canonical scalar product} (or simply {\it scalar product}).
The basic properties of the scalar product are easily derived from
Proposition 4 (in particular, the scalar product is commutative
if the $b$-semialgebra $A$ is commutative). The following theorem is the
idempotent analog of the Riesz--Fisher theorem.

{\bf Theorem 3.} {\sl Let a $b$-semialgebra $A$ be a semifield.
Then any nonzero $a$-linear functional $f$ on $A$ can be represented
as $f(y) = \langle y, x\rangle$, where $x\in A$, $x\neq \0$ and
$\langle \cdot, \cdot\rangle$ is the canonical scalar product on} $A$.

{\bf Remark.} Using the completion precedures, one can extend all
the results obtained to the case of incomplete semirings, spaces,
 and semimodules.

{\bf 6.} {\bf Example}. Let ${\cal B} (X)$ be a set of all bounded
functions with values belonging to $\bR(\max, +)$ on an arbitrary 
set $X$ and let $\widehat{\cal B}(X) = {\cal B} (X)\cup \{\0\}$.  
The pointwise idempotent addition of functions
$(\varphi_1 \oplus \varphi_2) (x) = \varphi_1(x)\oplus \varphi_2(x)$
and the multiplication $(\varphi_1\odot \varphi_2)(x) = (\varphi_1(x))\odot
( \varphi_2(x))$ define on $\widehat{\cal B} (X)$ the structure of a
$b$-semialgebra over the $b$-complete semifield $\bR_{\max}$. In
this case, $\1^* (\varphi) =\sup_{x\in X} \varphi (x)$  
and the scalar product is expressed in terms of idempotent integration
(see [1--3]): $\langle\varphi_1, \varphi_2\rangle  = 
\sup_{x\in X} (\varphi_1 (x)\odot\varphi_2 (x)) =  
\sup_{x\in X} (\varphi_1 (x) + \varphi_2 (x))
=\int\limits^{\oplus}_X(\varphi_1 (x)\odot\varphi_2 (x))\; dx$.
Scalar products of this type were systematically used in
idempotent snslysis. Using Theorems 1 and 3, one can easily
describe $a$-linear functionals on idempotent spaces in terms of
idempotent measures and integrals.

\medskip

\begin{center}
ACKNOWLEDGMENTS
\end{center}

This work was supported by the INTAS and Russian Foundation for
Basic Research, joint grant no. 95--91.

\medskip

\begin{center}
REFERENCES
\end{center}

 1. {\it Idempotent analysis}, Adv. Sov. Math., Maslov, V.P. and
 Samborski\u\i, \ S.N., Eds., Providence: Am. Math. Soc., 1992, vol. 13.

2. Maslov, V.P. and Kolokol'tsov, V.N., {\it Idempotentnyi analiz i ego
primenenie v optimal'nom upravlenii} (Idempotent Analysis and Its 
Application in Optimal Control), Moscow: Nauka, 1994.

3.  Kolokoltsov, V.N. and Maslov, V.P., {\it Idempotent Analysis and 
Applications}, Dordrecht: Kluwer, 1997.

4. Bacelli, F.L., Cohen, G., Olsder, G.J., and Quadrat, J.-P.,
{\it Synchronization and Linearity: An Algebra for Discrete Event Systems},
New York: Wiley, 1992.

5. Birkhoff, G., {\it Lattice Theory}, 
Amer. Math. Soc. Colloq. Publ., Providence, 1948. Translated under the title
{\it Teoriya reshetok}, Moscow: Nauka, 1984.

6. Shaefer, H., {\it Topological Vector Spaces}, New York:
Macmillan, 1966. Translated under the title {\it Topologicheskie
vektornye prostranstva}, Moscow: Mir, 1971.

\end{document}